\documentclass{amsart}

\usepackage{amsmath} 
\usepackage{amssymb}

\newtheorem{theorem}{Theorem}

\theoremstyle{definition}

\theoremstyle{remark}

\newtheorem{gq}[theorem]{Grinblat's Question}

\numberwithin{equation}{section}
\newcommand{\forces}{\Vdash}

\newcommand{\bV}{{\bf V}} 
\newcommand{\lesdot}{\mathrel{\mathord{<}\!\!\raise 
0.8 pt\hbox{$\scriptstyle\circ$}}} 
  
\newcommand{\cov}{{\rm {\bf cov}}\/}

\newcommand{\NULL}{{\rm Null}}
\newcommand{\meagre}{{\rm Meager}}

\newcommand{\bP}{{\mathbb P}}
\newcommand{\bQ}{{\mathbb Q}}
\newcommand{\bR}{{\mathbb R}}

\newcount\skewfactor
\def\mathunderaccent#1#2 {\let\theaccent#1\skewfactor#2
\mathpalette\putaccentunder}
\def\putaccentunder#1#2{\oalign{$#1#2$\crcr\hidewidth
\vbox to.2ex{\hbox{$#1\skew\skewfactor\theaccent{}$}\vss}\hidewidth}}
\def\name{\mathunderaccent\tilde-3 }

\begin{document}

\title[Consistency of ...]{Consistency of ``the ideal of null
restricted to some $A$ is $\kappa$--complete not $\kappa^+$--complete,
$\kappa$ weakly inaccessible and $\cov({\rm meagre})=\aleph_1$''}

\author{Saharon Shelah}
\address{Institute of Mathematics\\
 The Hebrew University of Jerusalem\\
 91904 Jerusalem, Israel\\
 and  Department of Mathematics\\
 Rutgers University\\
 New Brunswick, NJ 08854, USA}
\email{shelah@math.huji.ac.il}
\urladdr{http://www.math.rutgers.edu/$\sim$shelah}
\thanks{The research of the second author was partially supported by the
 Israel Science Foundation. Publication E52} 

\maketitle

In this note we give an answer to the following question of Grinblat (Moti
Gitik asked about it in the Oberwolfach meeting:

\begin{gq}
Is it consistent that 
\begin{enumerate}
\item[$(**)$] for some set $X$, $\cov(\NULL\restriction X)=\lambda$ is a
weakly inaccessible cardinal (so $X$ not null of course) while
$\cov(\meagre)$ is small, say it is $\aleph_1$.
\end{enumerate}
\end{gq}

\noindent {\sc A. The Forcing:}

Starting with a universe $\bV$ and a cardinal $\lambda$ of cofinality
$>\aleph_0$, regular for simplicity (otherwise the only difference is that
$J$ consists of ``bounded subsets''), in fact weakly inaccessible for
Grinblat's question.

Let $\bP=\bP_\lambda$ be the result of FS iteration $\langle\bP_i,
\name{\bQ}_i:i<\lambda\rangle$ with $\name{\bQ}_{2i}$ being the random real 
forcing, and $\name{\bQ}_{2i+1}$ being the Cohen forcing notion. Let
$\name{\bR}$ be a $\bP$--name for the forcing notion adding $\aleph_1$
random reals (i.e., forcing with the measure algebra of Borel subsets of
${}^{\omega_1} 2$ of positive Lebesgue measure).

We claim that $\bV_2=\bV^{\bP*\name{\bR}}$ is as required.

Let $\bV_1=\bV^\bP$.

As the whole forcing satisfies the ccc, no cardinal is collapsed etc
\bigskip 

\noindent {\sc B. Why $\cov(\meagre)=\aleph_1$ ?}

As forcing by $\name{\bR}$ does it (well known).
\bigskip 

\noindent {\sc C.}

Let  $\name{\eta}_i$ be the $\name{\bQ}_i$--generic real for
$i<\lambda$. Clearly they are pairwise distinct. Let 
\[\name{X}\stackrel{\rm def}{=}\{\name{\eta}_{2i}:i<\lambda\}.\]
This is a set of cardinality $\lambda$. Let $J$ be the ideal of subsets of
$X$ of cardinality $<\lambda$ (it is a $\lambda$--complete ideal on $X$). 

It is enough to prove
\begin{enumerate}
\item[$(*)$] $J$ is equal to the ideal of null subsets of $X$.
\end{enumerate}
\bigskip 

\noindent {\sc C1.}

Now, for every $\alpha<\lambda$ the set $\{\name{\eta}_{2i}:i<\alpha\}$ is
null in $\bV_2$. Why? Because $\name{\eta}_{2\alpha+1}$ is Cohen over
$\bV^{\bP_ {2\alpha+1}}$ the universe to which the above set belongs and is
an inner model of $\bV_2$. 

This is enough to show that every member of $J$ is null.
\bigskip 

\noindent {\sc C2.}

For the other direction, let $\name{Y}$ be a $\bP*\name{\bR}$ name of an
unbounded subset of $\lambda$. We shall prove that 
\[\{\name{\eta}_{2i}: i\in \name{Y}\}\] 
is forced to be non-null (this clearly suffices).

Let $p$ be a condition in $\bP*\name{\bR}$ forcing the inverse, so for some
$\bP*\name{\bR}$--name $\name{Z}$ of a null Borel subset of ${}^{\omega} 2$,
we have 
\[p\forces\mbox{`` }\{\name{\eta}_{2i}:i\in \name{Y}\}\subseteq \name{Z}
\mbox{ ''.}\]
We can find $\alpha<\lambda$ such that, in $\bV^{\bP_\alpha}$, $\name{Z}$
becomes an $\bR^{\bV^{\bP_\alpha}}$--name and $p$ is a member of
$\bR^{\bV^{\bP_\alpha}}$. 

Now for every $i$, if $\alpha<2i<\lambda$ then $\name{\eta}_{2i}$ is random
over $\bV^{\bP_\alpha}$. Hence, by the Fubini theorem (i.e., random reals
commute), it is also random over
$(\bV^{\bP_\alpha})^{\bR^{\bV^{\bP_\alpha}}}$.  Consequently it does not
belongs to $\name{Z}$, so we are done.

\end{document}